\documentclass[preprint]{imsart}

\usepackage[bookmarks=false]{hyperref}
\usepackage{amssymb}
\usepackage{amsfonts}
\usepackage{lineno} 
\usepackage{float}
\usepackage{mathrsfs}
\usepackage{booktabs}

\startlocaldefs

\newtheorem{theorem}{Theorem}
\newtheorem{example}{Example}
\newtheorem{remark}{Remark}
\newcommand{\proof}{\noindent{\bf Proof.  }}
\newcommand{\eproof}{$\Box$}
\endlocaldefs

\begin{document}

\begin{frontmatter}

\title{Marcinkiewicz Law of Large Numbers for Outer-products of Heavy-tailed, Long-range Dependent Data}
\runtitle{MSLLN for outer-products of Linear Models}

\begin{aug}
\author{\fnms{Michael A.} \snm{Kouritzin}\ead[label=e1]{michaelk@ualberta.ca}}
\and
\author{\fnms{Samira} \snm{Sadeghi}\thanksref{t1}\ead[label=e2]{ssadeghi@ualberta.ca}}

\thankstext{t1}{Corresponding author. E-mail: ssadeghi@ualberta.ca}
\runauthor{M.A. Kouritzin and S. Sadeghi}

\affiliation{University of Alberta}

\address{Department of Mathematical and Statistical Sciences\\
University of Alberta, Edmonton, AB  T6G 2G1 Canada }

\end{aug}

\begin{abstract}
\ The Marcinkiewicz Strong Law, 
$\displaystyle\lim_{n\to\infty}\frac{1}{n^{\frac1p}}\sum_{k=1}^n (D_{k}- D)=0$ a.s.\ 
with $p\in(1,2)$, 
is studied for outer products
$D_k=X_k\overline{X}_k^T$, where $\{X_k\},\{\overline{X}_k\}$ are both two-sided (multivariate) linear processes ( 
with coefficient matrices $(C_l), (\overline{C}_l)$ and i.i.d.\ zero-mean innovations $\{\Xi\}$, $\{\overline{\Xi}\}$).
Matrix sequences $C_l$ and $\overline{C}_l$ can decay slowly enough (as $|l|\to\infty$) that 
$\{X_k,\overline{X}_k\}$ have long-range dependence while $\{D_k\}$ can have heavy tails.
In particular, the heavy-tail and long-range-dependence phenomena for $\{D_k\}$ are handled simultaneously and 
a new decoupling property is proved that shows the convergence rate 
is determined by the worst of the heavy-tails
or the long-range dependence, but not the combination.
The main result is applied to obtain Marcinkiewicz Strong Law of Large Numbers
for stochastic approximation, non-linear functions forms and autocovariances.
\end{abstract}

\begin{keyword}[class=MSC]
\kwd[Primary ]{62J10}
\kwd{62J12}
\kwd{60F15}
\kwd[; secondary ]{62L20}
\end{keyword}

\begin{keyword}
\kwd{covariance}
\kwd{linear process}
\kwd{Marcinkiewicz strong law of large numbers}
\kwd{heavy tails}
\kwd{long-range dependence}
\kwd{stochastic approximation}
\end{keyword}

\end{frontmatter}

\section{Intoduction}
Let $D_k=X_k\overline{X}_k^T$ be random matrices with $\{X_k\}$, $\{\overline{X}_k\}$ 
being $\mathbb R^d$-valued (possibly two-sided, multivariate) linear processes 
\begin{eqnarray}\label{xylinearproc}
X_k = \sum_{l=-\infty}^\infty C_{k-l}\Xi_l, \quad \overline{X}_k = \sum_{l=-\infty}^\infty \overline{C}_{k-l}\overline{\Xi}_l.
\end{eqnarray}
defined on some probability space $(\Omega, F,P)$. 
$$\left\{\left(\Xi_{l}=( \xi_l^{(1)}, . . ., \xi_l^{(m)}), \overline{\Xi}_{l}=( \overline{\xi}_l^{(1)}, . . ., \overline{\xi}_l^{(m)})\right),
\ l\in{\mathbb Z}\right\}$$ 
are i.i.d.\ zero-mean random $\mathbb R^{m+m}$-vectors (innovations) such that 
$E[|\Xi_1|^2]<\infty$, $E[|\overline{\Xi}_1|^2]<\infty$ 
and $(C_{l})_{l\in\mathbb Z}$, $(\overline{C}_{l})_{l\in\mathbb Z}$ are $\mathbb R^{d\times m}$-matrix sequences
 satisfying 
$\sup\limits_{l\in\mathbb Z}|l|^\sigma\|C_l\|<\infty$, $\sup\limits_{l\in\mathbb Z}|l|^{\overline{\sigma}}\|\overline{C}_l\|<\infty$ 
for some $ (\sigma,\overline{\sigma})\in \left(\frac{1}{2},1\right]$.
Hence, $\{D_k\}$ can have heavy tails as well as long-range dependence.  

Linear process models are heavily used in finance, engineering,
econometrics, and statistics. 
In fact, classical time-series theory mainly 
involves the statistical analysis of stationary linear processes. 
Current applications in network theory and financial mathematics
leads us to study time series models where $\{D_k\}$ can have heavy tails and long memory. 
Heavy-tailed data exhibits frequent
extremes and infinite variance, while positively-correlated long memory data displays great serial 
momentum or inertia.  
Heavy-tailed data with long-range dependence has been observed in a plethora of empirical data set over the last fifty years and so.  
For instance, Mandelbrot \cite{Mandel68} observed that long memory time series often were
heavy-tailed and self-similar. 

The possible rates of the convergence is affected by both long-range dependence and heavy-tailed.
There are two broad types of dependence for linear processes. 
If the coefficients $(C_l)$ are absolutely summable and innovations have second moments, then the covariances of $X_k$ are
summable and we say that $\{X_k\}$ is short-range dependence (SRD). 
On the contrary, we generically say that $\{X_k\}$ is long-range
dependence (LRD) if its covariances are not absolutely summable. Practically, 
by choosing appropriate coefficients, 
matrix sequence $(C_l)$ can decay slowly enough (as $|l|\to\infty$) such that 
$\{X_k\}$ shows LRD. We consider $\{D_k\}$ to have LRD too in this $\{C_l\}$ non-summable case even though the second moments for $D_k$ may not exist.
There are also two general kinds of randomness.
If each $D_k$ fails to have a second moment, then we say
it has heavy-tailed (HT) and is otherwise light-tailed (LT).
In our setting, $D_k$ will either have HT or LT depending upon the moments
of and dependence between $\Xi_1$ and $\overline \Xi_1$.

There few general Marcinkiewicz Strong Law of Large Numbers (MSLLN) results for partial sums of $X_k$ under both heavy-tailed
and the long-range dependence and the MSLLN for partial sums of nonlinear functions of $X_k$ 
is almost untouched. 
Our purpose here is to establish a method and a structure under which certain MSLLN for heavy-tailed
and the long-range-dependent phenomena can be handled properly. 
Technically, our goal is to prove:  
\[
\lim\limits_{n\rightarrow \infty}\frac{1}{n^{\frac{1}{p}}}\sum\limits_{k=1}^{n}\left(D_{k}-D\right)=0\ \quad \mbox{a.s.} \quad \mbox{for}\quad p<\frac1{2-\sigma-\overline{\sigma}}\wedge\alpha\wedge2,
\]
\\
when $\max\limits_{1\le i,j\le m}\sup\limits_{t\ge 0}t^{\alpha}P(|\xi_1^{(i)}\overline{\xi}^{(j)}_1|>t)<\infty$ for some $\alpha>1$ and 
$\sup\limits_{l\in\mathbb Z}|l|^\sigma\|C_l\|<\infty$, $\sup\limits_{l\in\mathbb Z}|l|^{\overline{\sigma}}\|\overline{C}_l\|<\infty$ when $ (\sigma,\overline{\sigma})\in \left(\frac{1}{2},1\right]$. 
This format of $\{D_k\}$ is critical for our result since, it allows  LRD and HT conditions decouple and convergence rate 
be determined by the worst of the HT requirement $p<(\alpha\wedge 2)$ and
the LRD condition $p<\frac1{2-\sigma-\overline{\sigma}}$, 
but not the combination. 
A bifurcation happens. Consider the summation, $\displaystyle D_k=\sum_{l,m=-\infty}^\infty C_{k-l}\Xi_l\overline{C}_{k-m}\overline{\Xi}_m$, 
broken into off-diagonal and diagonal terms. 
Due to the independence of $(\Xi_l,\overline \Xi_l)$ from
$(\Xi_m,\overline \Xi_m)$, the off-diagonal sum
\(
\sum\limits_{l\neq m}C_{k-l}\overline C_{k-m}\Xi_l\overline \Xi_m 
\)
does not have heavy tails ( when  $\alpha>1$ ).
Conversely, since $\sigma+\overline \sigma>1$ the diagonal sum 
\(
\sum\limits_{l=-\infty}^\infty C_{k-l}\overline C_{k-l}\Xi_l\overline \Xi_l 
\)
does not experience long-range dependence. In addition, the rate of convergence depends on the worst of $(\alpha\wedge 2)$ and $\frac1{2-\sigma-\overline{\sigma}}$, 
so whenever we are in the LRD dominant case, ($\alpha>\frac1{2-\sigma-\overline{\sigma}}$), the off-diagonal terms dictate the rate of convergence by the LRD effect ($p<\frac1{2-\sigma-\overline{\sigma}}$) and in the 
HT dominant case, ($\alpha<\frac1{2-\sigma-\overline{\sigma}}$), the diagonal terms dictate the rate of convergence by HT effect ($p<\alpha$). The bifurcation point is when $\alpha=\frac1{2-\sigma-\overline{\sigma}}$ and $\alpha<2$.  
 
\section{Background} 
In this section we give a review of some existing literature on MSLLN or weak convergence for 
partial sums, sample covariance and non-linear function of partial sums with heavy-tailed and/or long-range dependence. 
Many existing results were only established in the scalar case. 
For ease of assimilation we use $\{x_k\}$, $(c_l)$, $\{d_k\}$ and $\{\xi_k\}$ to denote these scalar versions of $\{X_k\}$, $(C_l)$, $\{D_k\}$ and $\{\Xi_k\}$ and 
$\{x_{k+h}\}$ for $\{\overline{X}_k\}$ when it is a shifted version of $\{x_k\}$.

\subsection{Partial Sums}
There are many of publications that consider almost sure rates of convergence for linear processes 
having either LRD or HT. 
However, there are only a few like Louhchi and Soulier \cite{30} that 
considered the combination of these two phenomena.
They stated the following result for
linear symmetric $\alpha$-stable (S$\alpha$S) processes. 
\begin{theorem}\label{linear}
Let $\{\xi_j\}_{j\in\mathbb Z}$ be i.i.d. sequence of S$\alpha$S random variables
with $1<\alpha<2$ and $\{c_j\}_{j\in\mathbb Z}$ be a bounded collection such that 
\(
\sum\limits_{j\in\mathbb Z}|c_j|^s<\infty
\)
for some $s\in[1,\alpha)$. 
Set
\(
x_k=\sum\limits_{j\in\mathbb Z}c_{k-j}\xi_j.
\)
Then, 
for $p\in(1,2)$ satisfying $\frac1p>1-\frac1s+\frac1\alpha$ 
\[
\frac1{n^\frac1p}\sum_{i=1}^n x_i\to 0\ \mbox{a.s.}
\]
\end{theorem}
The condition $s<\alpha$ ensures 
\(
\sum\limits_{j\in\mathbb Z}|c_j|^\alpha<\infty
\)
and thereby convergence of 
\(
\sum\limits_{j\in\mathbb Z}c_{k-j}\xi_j.
\)
Moreover, $\{x_k\}$ not only exhibits heavy tails but also long-range dependence
if, for example, $c_j=|j|^{-\sigma}$ for $j\neq 0$ and some $\sigma\in \left(\frac12,1\right)$.
Notice there is interactions between the heavy tail condition and the long range dependent 
condition.
In particular for a given $p$, heavier tails ($\alpha$ becomes smaller) 
implies that you can not have as 
long range dependence ($s$ must become smaller) and vice versa. 
Moreover, this result is difficult or even impossible to apply  in our outer product setting  due to the fact that $x_k$'s are linear processes with  
S$\alpha$S innovations and so $x_k$ cannot be decomposed to product of two variables even in the scalar case.

\subsection{Non-linear function of partial sums}
The limit behavior of suitably normalized partial sums of stationary random variables that demonstrate either LRD
 or HT has been subject of study by many authors.  
Applications can be found in geophysics, economics, hydrology and statistics.
For instance, in contexts like Whittle approximation, the asymptotic behavior of quadratic forms of stationary sequences 
have an important role. In addition, the efficacy of ``$R/S$-statistic'' theory that was introduced 
for estimating the long-run, non-periodic statistical dependence of time series
by Hurst and developed by Mandelbrot \cite{Mandel72}, can be confirmed by convergence of these limit functions. 

There are many results that deal with  the existence and description of limit distributions of sums 
\begin{eqnarray}\label{SN}
S_{n,h}(t)=\sum\limits _{k=1 }^{[nt] } (h(x_k)-E(h(x_k))), \quad t\geq 0,
\end{eqnarray}
where $h$ is a (nonlinear) function.  
The limit behavior for a Gaussian LRD process $\{x_k\}$, firstly was studied by Rosenblatt
 \cite{Rosen}. Afterward, Dobrushin and Major \cite{DM} explained it in more general form. 
Then Taqqu \cite{Taqqu79} showed that the limit in distribution of particular normalized
sums $S_{n,h}(t)$ is determined by the Hermite rank $m^*\in\{1,2,...\}$ of $h(x)$,  
which is the index of the first nonzero coefficient in the Hermite expansion. 
On the other hand, the behavior of nonlinear non-Gaussian LRD processes is much
less commonly known. One of the most studied models of non-Gaussian LRD processes
is the one-sided linear (moving average) process, 
\begin{eqnarray}\label{X}
x_k=\sum\limits _{j=0 }^{\infty} c_j \xi_{k-j}, 
\end{eqnarray}
 in which, innovations $\xi_k, k\in \mathbb Z$, are independent and identically distributed
(i.i.d.), have zero mean with finite variance,   
and coefficients $c_j $ satisfy:
\begin{eqnarray}\label{b}
c_j \sim c_\sigma j^{-\sigma},  \quad\,\, j\geq1\,\, 
\end{eqnarray}  
for some constant $c_\sigma \neq 0$, $ c_0 = 1 $ and $\sigma\in (\frac12, 1)$.\\

Surgailis \cite{SU82} considered the limit behavior of  partial sum processes
$S_{n,h}(t)$ of polynomial $h$ of linear process $\{x_k\}_{k\in \mathbb Z}$. Later,   
Giraitis and Surgailis \cite{GS86}\cite{GS89}, Avram and Taqqu  \cite{AT87} noticed that 
the only difference between this case and Gaussian case is that
the Hermite rank $m^*$ of $h(x)$ has to be replaced by the Appell rank $m$.

Vaiciulis \cite{vai2003} investigated distributional convergence for normalized
partial sums of Appell polynomials $A_m(x_k)$ of linear processes
$x_k$ having both long-memory and heavy-tails in the sense
$EA_m^2(x_k)=\infty$. In particular, he assumed $x_k$ had the form (\ref{X}) with innovations
$\{\xi_k^m\}$ belonging to the domain of attraction of an $\alpha$-stable law with
$1<\alpha<2$ and $c_j$ following (\ref{b}). The limit was: {\bf i)} an $\alpha$-stable Levy process, {\bf ii)} an $m^{th}$ 
order Hermite process, or {\bf iii)} the sum of two mutually independent
$\alpha$-stable Levy and $m^{th}$ order Hermite processes, depending on the value of $\alpha, m$ and $\sigma$ where $\sigma \in (\frac12, 1)$.

Thereafter, Surgailis \cite{S2004}  
considered the bounded, infinitely differentiable $h$ case where $\{x_k\}$ was LRD and had innovations with probability tail decay of $x^{-2\alpha}$ for $1<\alpha<2$.  
Suppose $x_k $ satisfies (\ref{X}) and (\ref{b}). Then he showed three different limiting behaviors corresponding to three 
different LRD-HT setting: 
 $n^{1-(2\sigma-1)m^*/2}S_{n,h}(t)$,  $n^{\frac{1}{2\alpha\sigma}}S_{n,h}(t)$ or $n^{\frac12}S_{n,h}(t)$ converge in distribution to respectively a  Hermite process of order $m^*$,
 a $2\alpha\sigma$-stable Levy process or a Brownian motion, all at time t, for certain range of $\alpha$ and $\sigma$.  

\subsection{Sample Covariances}
Auto-covariance functions play a substantial role in time series analysis and have diverse applications in 
inference problems, including hypothesis testing and parameter estimation. The natural estimator of auto-covariance
is sample covariance. Hence, the convergence properties of the sample covariance   
is of great interest. 
In the case of LRD and HT, it is an area of active research.

Davis and Resnick \cite{DR86}  studied the distributional convergence
 of sample autocovariances for two-sided linear processes with innovations that were
i.i.d. and had regularly varying tail probabilities of index
$\alpha >0$. 
\begin{eqnarray}\label{tdecay} \nonumber
&& P(|\xi_k|>x)=x^{-2\alpha} L(x),\\
&& \frac{P(\xi_k>x)}{P(|\xi_k|>x)}\to p \quad \mbox{and}\quad \frac{P(\xi_k<-x)}{P(|\xi_k|>x)}\to q, \quad \mbox{as}\,\, x\to \infty,
\end{eqnarray}
where $L(.)$ is a function slowly varying at infinity $\left(\displaystyle \lim_{j\to \infty} \frac{L(aj)}{L(j)}=1\right)$ and $0\leq p\leq1,\,\, q=1-p$.
They considered the case where the innovations had finite variance ($\iota$) but infinite
fourth moment, i.e. $ 1<\alpha <2$ with absolutely summable coefficients $c_j$ with form of (\ref{b}).

{\bf Note:} We choose to scale our constants, here and in the sequel, so that $\alpha<2$ always mean HT of the object of interest, which is $x_kx_{k+h}$ 
or more generally $X_k\overline{X}_k$. 

In case of infinite fourth moment for $\{\xi_k\}_{k\in \mathbb Z}$, the asymptotic distribution of normalized
sample autocovariances of long-memory processes was studied by  Horv{\'a}th and Kokoszka \cite{HK}.    
Suppose we observe the realization $x_1, x_2, ..., x_{n+v}, \, n>1 ,v\geq0$, the sample autocovariances and 
population autocovariances are defined as
\begin{eqnarray}\label{autocov}\nonumber
&&\hat{\gamma}_h^{(n)}=\frac1n \sum_{k=1}^{n} x_k x_{k+h}, \quad h=0,1,..., v, \quad\mbox{and}\\
&& \gamma_h= E[x_0x_h]=\iota \sum_{j=0}^{\infty} c_j c_{j+h},
\end{eqnarray}
respectively.
Horv{\'a}th and Kokoszka \cite[Theorem $3.1$]{HK} studied the asymptotic distribution $[\hat{\gamma}_h^{(n)}-\gamma_h], \, h=0,1,...,v$ 
for linear process of form (\ref{X}) with coefficients and innovations satisfying (\ref{b}) and (\ref{tdecay}) 
and a norming constant  $a_n=\inf \{x: P(|\xi_1|> x)\leq n^{-1}\}$ (roughly of order $n^\frac1{2\alpha}$) satisfying
\begin{eqnarray}\label{ct} 
\lim_{n\to \infty}nP[|\xi_k|>a_nx]=x^{-2\alpha},\,\, x>0.
\end{eqnarray} 
We quote this result in our notations as the following theorem.

\begin{theorem}\label{HTLRD}
Suppose, conditions (\ref{X}), (\ref{b}), (\ref{tdecay}) and (\ref{ct}) hold.
\begin{itemize}
\item[\bf (a)] If $1-\frac1{2\alpha} <\sigma<1$ and $1<\alpha<2$, then\\

$na_n^{-2}[\hat{\gamma}_h^{(n)}-\gamma_h]\stackrel{d}{\rightarrow}\left(S-\frac{\alpha}{\alpha-1}\right)\left[\displaystyle\sum_{j=0}^{\infty} c_j c_{j+h}\right],\quad h=0,1,...,H.$\\

where $S$ is an $\alpha$-stable random variable. For the above to hold for $\sigma=3/4$, we must additionally 
assume that $a_n^{-4}n\ln n\to 0$.\\

\item[\bf (b)] If $\frac12<\sigma<1-\frac1{2\alpha}$ and $1<\alpha<2$, then\\

$n^{2\sigma-1}[\hat{\gamma}_h^{(n)}-\gamma_h]\stackrel{d}{\rightarrow}\iota c_\sigma^2\left[U_\sigma(1)\right],\quad h=0.1,...,H.$\\

where $U_\sigma$ is a Rosenblatt process. The Rosenblatt process is often defined by the iterated stochastic integral:
\end{itemize}
$U_\sigma(t)=2 \int_{w_1<w_2<t}\left[ \int_{0}^{t}(\tau-w_1)_+^{-\sigma}(\tau-w_2)_+^{-\sigma}d\tau\right]W(dw_1)W(dw_2)$,\\ 
in which $W(.)$ is the standard Wiener process on the real line.
\end{theorem}

This theorem works for one-sided linear processes with a regularly varying tail condition and  
gives us weak convergence. 

Notice that in Theorem \ref{HTLRD}, case {\bf (a)} represents the HT dominant, ($\alpha<\frac1{2-2\sigma}$), 
so the diagonal terms dictate convergence to an $\alpha$-stable distribution. However,  
case {\bf (b)}  represents the LRD dominant, ($\alpha>\frac1{2-2\sigma}$), hence off-diagonal terms take over and we get convergence to Rosenblatt process.

\section{Main results}

Our first result is in the scalar case.
Later, we will extract the full vector-valued result as a second main theorem.
All proofs are delayed until the next section after we have discussed the applications.

\begin{theorem}\label{HeavyLong}
Let $ \left\{(\xi_{l},\overline \xi_{l})\right\}_{l\in \mathbb Z}$ be i.i.d.\ zero-mean random variables such that  
$E[\xi_1^2]<\infty$, $E[\overline \xi_1^2]<\infty$ and
$\sup\limits_{t\ge 0 }t^{\alpha}P(|\xi_1\overline \xi_1|>t)<\infty$ 
for some $\alpha>1$.
Moreover, suppose 
$(c_{l})_{l\in\mathbb Z}, (\overline c_{l})_{l\in\mathbb Z}$ satisfy
\begin{eqnarray*}
\sup\limits_{l\in\mathbb Z}|l|^\sigma|c_l|<\infty,\quad 
\sup\limits_{l\in\mathbb Z}|l|^{\overline \sigma}|\overline c_l|<\infty \quad
\mbox{for some}\quad \sigma, \overline \sigma \in \left(\frac{1}{2},1\right], 
\end{eqnarray*}
$d_k=\sum\limits_{l,m=-\infty}^{\infty}c_{k-l}\overline c_{k-m}\xi_{l}\overline \xi_{m}$ 
and $d=E[\xi_1\overline \xi_1]\sum\limits_{l=-\infty }^{\infty }c_{l}\overline c_{l}$. 
Then, for $p$ satisfying $p<\frac{1}{2-\sigma -\overline \sigma}\wedge \alpha\wedge 2$ 
\[
\lim\limits_{n\rightarrow \infty}\frac{1}{n^{\frac{1}{p}}}\sum\limits_{k=1}^{n}\left(d_{k}-d\right)=0\ \mbox{a.s.} 
\]
\end{theorem}
\begin{remark}
The tail probability bound ensures that $E[|\xi_1\overline \xi_1|^r]<\infty$ for
any $r\in(1,(\alpha\wedge2) )$ and $E[d_1]$ exists but 
it is possible that $E[d_1^2]=\infty$ so we are handling heavy tails for $\{d_k\}$.
On the other hand, $E[|\xi_1\overline \xi_1|^\alpha]<\infty$ implies our tail condition
by Markov's inequality.
$\sigma$, $\overline\sigma$ bound the amount of long-range dependence in 
\(
x_{k}=\sum\limits _{l=-\infty }^{\infty }c_{k-l}\xi _{l}\),
\(
\overline x_{k}=\sum\limits _{l=-\infty }^{\infty }\overline c_{k-l}\overline \xi _{l}\).
If $\sigma$ can be taken larger than $1$, then 
\(
\sum\limits_{k=1}^\infty E[x_0x_k]<\infty
\)
and there is no long-range dependence in $\{x_k\}$.
$\sigma> \frac12$ with $E[\xi_1^2]<\infty$ ensures that 
$\sum\limits _{l=-\infty }^{\infty }c_{k-l}\xi _{l}$ converges a.s.
\end{remark}

\begin{remark}
Notice that the constraints to handle long-range dependence, $p<\frac{1}{2-\sigma -\overline \sigma}$, and to 
handle the heavy tails, $p<(\alpha\wedge2)$, \emph{decouple}.
This decoupling appears to be due to the structure of $d_k$.
Due to the independence of $(\xi_l,\overline \xi_l)$ from
$(\xi_m,\overline \xi_m)$, the off-diagonal sum
\(
\sum\limits_{l\neq m}c_{k-l}\overline c_{k-m}\xi_l\overline \xi_m 
\)
does not have heavy tails.
Conversely, since $\sigma+\overline \sigma>1$ the diagonal sum 
\(
\sum\limits_{l=-\infty}^\infty c_{k-l}\overline c_{k-l}\xi_l\overline \xi_l 
\)
does not experience long-range dependence.
\end{remark}

We will give a simple example to verify conditions in Theorem \ref{HeavyLong}.
Recall, a non-negative random variable $\xi$ obeys a power law with parameters 
$\beta>1$ and $x_{\min}>0$, written $\xi\sim PL (x_{min},\beta)$, if it has density 
\begin{eqnarray*}
f(x)=\frac{\beta-1}{x_{min}} (\frac{x}{x_{min}})^{-\beta} \quad \forall \,\, x\geq x_{min}
\end{eqnarray*}
so $E|\xi|^r= \left\{\begin{array}{ll} x_{min}^r (\frac{\beta-1}{\beta-1-r})& \quad r<\beta-1\\
\infty & \quad r\ge\beta-1\end{array}\right.$.\\
It has a folded $t$ distribution with 
parameter $\beta>1$, written $\xi\sim Ft (\beta)$, if it has density  
\begin{eqnarray*}
f(x)=\frac{2\Gamma( \frac{\beta}{2})}{\Gamma( \frac{\beta-1}{2})\sqrt{(\beta-1) \pi}} \left( 1+\frac{x^2}{(\beta-1)} \right) ^{-\frac{\beta}{2}} \quad \forall \,\, x> 0
\end{eqnarray*}
so $E(|\xi|^r)$ exists if and only if $r<\beta-1$. 
\begin{example}
Suppose $p,q,\alpha,\beta,\overline \beta>1$ are such that $\frac1p+\frac1q=1$, $\beta>p\alpha+1$,
$\overline \beta>q\alpha+1$
and $p\alpha,q\alpha\ge 2$. 
If $\xi_1$ and $\overline{\xi_1}$ have power law distribution, lets say $\xi_1\sim Pl (x_{\min},\beta)$, $\overline \xi_1\sim Pl (\overline x_{\min},\overline \beta)$
for some $x_{\min},\overline x_{\min}>0$, then 
$E[\xi_1^2],\ E[\overline \xi_1^2]<\infty$ and 
$\sup\limits_{t\ge 0 }t^{\alpha}P(|\xi_1\overline \xi_1|>t)<\infty$. 
If $\xi_1\sim Ft (\beta)$, $\overline \xi_1\sim Ft (\overline \beta)$, then 
$E[\xi_1^2],\ E[\overline \xi_1^2]<\infty$ and 
$\sup\limits_{t\ge 0 }t^{\alpha}P(|\xi_1\overline \xi_1|>t)<\infty$. 
Either way, the Theorem \ref{HeavyLong} applies with properly chosen $(c_l,\overline{c}_l)$.
\end{example}

We now consider the case where $X_k$ and $\overline{X}_k$ are (multivariate) linear processes.

\begin{theorem}\label{MatrixHeavyLong}
Let $ \left\{\Xi_{l}\right\}$ and $ \left\{\overline{\Xi}_{l}\right\}$ be i.i.d.\ zero-mean random $\mathbb R^m$-vectors such that 
$\Xi_{l}=\left( \xi_l^{(1)}, . . ., \xi_l^{(m)}\right) $, $\overline{\Xi}_l=\left( \overline{\xi}_l^{(1)}, . . ., \overline{\xi}_l^{(m)}\right) $, 
 $\max \limits_{1\le i,j\le m} \sup\limits_{t\ge 0}t^{\alpha}P(|\xi_1^{(i)}\overline{\xi}^{(j)}_1|>t)<\infty$
for some $\alpha >1$, $E[|\Xi_1|^2]<\infty$ and  $E[|\overline{\Xi}_1|^2]<\infty$.
Moreover, suppose matrix sequences
$(C_{l})_{l\in\mathbb Z}, (\overline{C}_{l})_{l\in\mathbb Z}$ $\in \mathbb R^{d\times m}$ satisfy
\begin{eqnarray*}
\sup\limits_{l\in\mathbb Z}|l|^\sigma\|C_l\|<\infty,\, \sup\limits_{l\in\mathbb Z}|l|^{\overline{\sigma}}\|\overline{C}_l\|<\infty \quad 
\mbox{for some} \quad (\sigma,\overline{\sigma}) \in \left(\frac{1}{2},1\right], 
\end{eqnarray*}
$X_k$, $\overline{X}_k$ take form of (\ref{xylinearproc}), $D_k=X_k\overline{X}_k^T$ and $D=E[X_1\overline{X}_1^T]$.
Then, for $p$ satisfying $p<\frac{1}{2-\sigma -\overline \sigma }\wedge\alpha\wedge2$ 
\[
\lim\limits_{n\rightarrow \infty}\frac{1}{n^{\frac{1}{p}}}\sum\limits_{k=1}^{n}\left(D_{k}-D\right)=0\ \quad \mbox{a.s.} 
\] 
\end{theorem}
This theorem follows by linearity of limits and Theorem \ref{HeavyLong}.

\subsection{Applications}
We give some applications of our theorems. 

\subsubsection{Stochastic Approximation}
Stochastic approximation (SA) is often used in optimization problems for linear models.
Hence, the convergence properties of SA algorithms driven by linear models is of
utmost interest.
For illustration, we assume $\{z_k, k=1,2,..\}$ and $\{y_{k},k=2,3,...\}$ are respectively
$\mathbb R^d-$ and $\mathbb R-$valued stochastic processes, defined on some probability 
space $(\Omega, F,P)$, that satisfy 
\begin{eqnarray}\label{observeqn}
y_{k+1}=z_k^Th+\epsilon_k,\quad\quad \forall k=1,2,\dots,
\end{eqnarray}
where $h$ is an unknown $d$-dimensional parameter or weight vector of interest and $\{\epsilon_k\}$ is a noise sequence.
We want to estimate the parameter vector $h$ through the stochastic approximation algorithm:
\begin{eqnarray}\label{a3}
	h_{k+1} = h_k +\mu_k (b_k-A_kh_k), 
\end{eqnarray}
where $\mu_k$ is the $k^{\rm th}$ step gain of the form 
$\mu_k=k^{-\chi}$ for some $\chi\in\left(\frac12,1\right]$, $A_k=z_kz_k^T$ and $b_k=y_{k+1}z_k$.\\
Kouritzin and Sadeghi \cite{MS} studied the convergence and almost sure rates of convergence for the  algorithm (\ref{a3}).  
Now, we can combine our main result (Theorem \ref{MatrixHeavyLong} ) with \cite[Corollary $2$]{MS} to obtain a powerful 
rate of convergence result for stochastic approximation.
\begin{theorem}\label{LinearMain}
Let $ \left\{\Xi_{l}\right\}$ be i.i.d.\ zero-mean random $\mathbb R^m$-vectors such that 
$$\sup\limits_{t\ge 0 }t^{\alpha}P(|\Xi_1|^2>t)<\infty\quad \mbox{for some}\,\, \alpha\in(1,2)$$
$(C_{l})_{l\in\mathbb Z}$ be $\mathbb R^{(d+1)\times m}$-matrices such that
$\sup\limits_{l\in\mathbb Z}|l|^\sigma\|C_l\|<\infty\, \mbox{for some}\, \sigma\in \left(\frac{1}{2},1\right]$,  
\begin{eqnarray}\nonumber
(z_k^T,y_{k+1})^T = \sum_{l=-\infty}^\infty C_{k-l}\Xi_l,
\end{eqnarray}
$A_k=z_kz_k^T$ and $b_k=y_{k+1}z_k$ and
$A=E[z_kz^T_k]$ and $b=E[y_{k+1}z_k]$.\\
Then, $|h_n- h|=o(n^{-\gamma})$ as $n\to \infty$ a.s. for any 
$\gamma<\gamma_0^{(\chi)}\doteq(\chi-\frac1\alpha)\wedge(\chi+2\sigma-2)$.
\end{theorem}

\proof By Theorem \ref{MatrixHeavyLong} when $\frac1p=\chi-\gamma$, 
$\overline{X}_k^T=X_k^T=(z_k^T,y_{k+1})$, $\overline{\Xi}_{l}=\Xi_l$, 
$\overline{C}_{l}=C_l$, $\overline{\sigma}=\sigma$, and
\(
D_k=\left(\begin{array}{l c r}
	z_kz_k^T&&y_{k+1}z_k\\
	y_{k+1}z_k^T&&y_{k+1}^2
\end{array}\right), 
\)
$$\frac{1}{n^{\chi-\gamma}}\sum\limits_{k=1}^n (D_{k}- D)\to 0\quad a.s.,$$
 where 
\(
D=\left(\begin{array}{l c r}
	A&&b\\
	b^T&&E[y_{k+1}^2]
\end{array}\right).  
\)
 The first $d$-rows of $\frac{1}{n^{\chi-\gamma}}\sum\limits_{k=1}^{n}\left(D_{k}-D\right)\to 0\ \mbox{a.s.}$ then establish the MSLLN 
$$\frac{1}{n^{\chi-\gamma}}\sum\limits_{k=1}^n (A_{k}- A)\to 0
\quad \mbox{and}\quad \frac{1}{n^{\chi-\gamma}}\sum\limits_{k=1}^n (b_k-b)\to 0 \quad a.s. $$
Now, we apply \cite[Corollary $2$]{MS} to complete the proof. \eproof

\begin{remark}
Note that $\chi-\gamma$ satisfies the required conditions 
$\chi-\gamma>2-2\sigma$ and $\chi-\gamma>\frac1\alpha$ in Theorem \ref{MatrixHeavyLong}. 
Theorem \ref{LinearMain} also appears in  \cite[Theorem $7$]{MS}.
\end{remark}

\subsubsection{Non-linear Function of linear processes}
As mentioned in Background, Vaiciulis \cite{vai2003} showed the convergence of distributions of the partial sum processes 
with non-linear $h(x_k)$ in terms of convergence of Appell polynomials $A_m(x_k)$ of a long-memory
moving average process $\{x_k\}$ with i.i.d. innovations $\{\xi_k\}$ in the case where the variance $EA_m^2(x_k)=\infty$, 
and the distribution of $\xi_1^m$ belongs to the domain of attraction of an $\alpha$-stable law with $1<\alpha<2$. 

Practically, the simplest examples of functions $h(x)$ with a given Appell rank $m$ are Appell polynomials $h=A_m$  
relative to the marginal distribution $x_1$ of the linear process (\ref{X}). In case $m=2$ the Appell polynomial is 
 $A_2(x)=x^2-\mu_2$ where $\mu_2=Ex^2$.
Viaiciulis \cite[Theorems $1.1$ and $1.2$]{vai2003} proved that when $m(2\sigma-1)<1,$ $m\geq2$ and $\sigma \in (\frac12, 1)$ the limit distribution of partial
sums of $m^{th}$ Appell polynomial is either {\bf (i)} an $\alpha$-stable Levy process for $2-2\sigma<1+\frac2m(\frac1\alpha-1)$, 
or {\bf (ii)} an $m^{th}$ order Hermite process for $2-2\sigma>1+\frac2m(\frac1\alpha-1)$ or {\bf (iii)} 
the sum of two mutually independent processes depending on the value of $\alpha, m$ and $\sigma$, for  $2-2\sigma=1+\frac2m(\frac1\alpha-1)$.

Taking into account all his conditions ( when $t=1$ ) and transforming it to our case we write our complementary almost sure rate-of-convergence theorem.

\begin{theorem}\label{our-quad}
Suppose $A_2$ represents the Appell polynomials with rank $2$ relative to the marginal distribution $x_1$ of 
the linear process $\displaystyle x_k=\sum\limits _{j=0 }^{\infty} c_{k-j} \xi_j$,  
for  $p\in [1,\frac{1}{2-2\sigma}\wedge\alpha)$ when 
\begin{eqnarray}\label{tail}
\sup\limits_{t\ge 0 }t^{\alpha}P(\xi_1^2>t)<\infty \quad\mbox{for some}\quad \alpha\in(1,2),
\end{eqnarray}
\begin{eqnarray}\label{lrd}
\sup\limits_{l\in\mathbb Z}|l|^\sigma|c_l|<\infty \quad\mbox{for some}\quad \sigma\in \left(\frac{1}{2},1\right].
\end{eqnarray}
Then, 
\[
\lim\limits_{n\rightarrow \infty}\frac{1}{n^{\frac{1}{p}}}\sum\limits_{k=1}^{n}A_2(x_k)=0\ \mbox{a.s.} 
\]
\end{theorem}

One might wonder if we have obtained the best possible MSLLN. Indeed, we have. 
For example for $m=2$, Viaiciulis \cite{vai2003} shows convergence in distribution of $\frac{1}{n^{(2-2\sigma)\wedge\frac1\alpha}} \sum\limits_{k=1}^{n}A_2(x_k)$ to different non-trivial limits 
in cases $(2-2\sigma)>\frac1\alpha$ (LRD dominant) or $(2-2\sigma)<\frac1\alpha$ (HT dominant), respectively.
Therefore, $\displaystyle \frac{1}{n^{(2-2\sigma)\wedge\frac1\alpha}}\sum\limits_{k=1}^{n}A_2(x_k)$ cannot converge to zero almost surely. Theorem \ref{our-quad} gives MSLLN for Appell polynomials with rank $2$ 
or in other word gives the convergence and almost sure rates of convergence for partial sums of second Appell polynomial  
when $\frac1p>(2-2\sigma) \vee \frac1\alpha$. 
Our result is optimal in polynomial sense and we cannot do better than that in terms of MSLLN.  

\subsubsection{Autocovariances}
As mentioned in the background, autocovariance estimation under HT and LRD conditions is
an active area of research.
We will handle the asymptotic behavior 
of sample covariance function for processes with LRD, innovations of  infinite $4^{th}$ moment and finite variance $\iota$.
If we define the sample aurtocovariance and population autocovariance functions by $\hat{\gamma}^{(n)}(h)$ and $\gamma(h)$, as (\ref{autocov}), 
we have following almost sure result.

\begin{theorem}\label{our-auto}
Assume $\hat{\gamma}^{(n)}(h)$ and $\gamma(h)$, as (\ref{autocov}) in which  $x_k=\sum\limits _{j=0 }^{\infty} c_{k-j} \xi_j$ 
and satisfies (\ref{tail}) and (\ref{lrd}) with $E[\xi_1^2]=\iota$.  
Then for  $p$ satisfying $p<\frac{1}{2-2\sigma}\wedge\alpha\wedge2$ 
\begin{eqnarray}\label{auto-app}
n^{1-\frac1p}[\hat{\gamma}_h^{(n)}-\gamma_h] \to 0 \,\, \mbox{a.s.}
\end{eqnarray}
\end{theorem} 
\proof Note that in Theorem \ref{HeavyLong}, for case $\overline{\xi}_l=\xi_l$, $E[\xi_1^2]=\iota$,   
$\overline{c}_l=c_{l+h}$ and $\{c_{l}=0,\, \forall \,l<0\}$ we have
 $$d_k=\sum\limits_{l=-\infty}^{k}\sum\limits_{m=-\infty}^{k+h}c_{k-l} c_{k+h-m}\xi_l \xi_{m}\,\quad \mbox{and}\quad d=\iota \sum\limits_{l=0}^{\infty}c_{l}c_{l+h}.$$

Hence, 
\begin{eqnarray}
\!\!\!\!&& \quad \frac{1}{n^{\frac1p}}\!\sum\limits_{k=1}^{n}\left(d_{k}-d\right)
=\frac{1}{n^{\frac1p}}\sum\limits_{k=1}^{n}\!\!\left(\sum\limits_{l=-\infty}^{k}\sum\limits_{m=-\infty}^{k+h}c_{k-l}c_{k+h-m}\xi_l\xi_m - \iota \sum\limits_{l=0}^{\infty}c_{l}c_{l+h} \right).
\end{eqnarray}
On the other hand, (\ref{auto-app}) can be written as 
\begin{eqnarray}\label{tcov}\nonumber
&&n^{1-\frac1p}[\hat{\gamma}_h^{(n)}-\gamma_h]=\frac{1}{n^{\frac1p}}\sum\limits_{k=1}^{n}(x_kx_{k+h}-Ex_0x_h)\\
&&=\frac{1}{n^{\frac1p}}\sum\limits_{k=1}^{n}\left(\sum\limits_{l=-\infty}^{k}\sum\limits_{m=-\infty}^{k+h}c_{k-l}c_{k+h-m}\xi_l\xi_m - \iota \sum\limits_{l=0}^{\infty}c_{l}c_{l+h} \right).
\end{eqnarray}
So, the result follows.\eproof

As we saw, Theorem \ref{HTLRD} gives the convergence to the following non trivial limits for $\frac{2\alpha-1}{2\alpha} <\sigma<1$ and  $\frac12<\sigma<\frac{2\alpha-1}{2\alpha}$ when $1<\alpha<2$, 
\begin{eqnarray}\label{as-auto}\nonumber
&\mbox{{\bf (a)}}&\quad \frac{1}{a_n^2}\sum\limits_{k=1}^{n}(x_kx_{k+h}-Ex_0x_h)\stackrel{d}{\rightarrow}\left(S-\frac{\alpha}{\alpha-1}\right)\left[\sum_{l=0}^{\infty} c_l c_{l+h}\right],\\\nonumber
&\mbox{{\bf (b)}}&\quad \frac{1}{n^{2-2\sigma}}\sum\limits_{k=1}^{n}(x_kx_{k+h}-Ex_0x_h)\stackrel{d}{\rightarrow}\iota c_\sigma^2\left[U_\sigma(1)\right],
\end{eqnarray}
respectively, for $h=0,1,...,v$. 

It is clear that in the case of HT dominant, $\frac1\alpha>2-2\sigma$, we have almost sure convergence (Theorem \ref{our-auto}) when $\frac1p> \frac1\alpha$.  
When $\frac1p=\frac1\alpha$, we get into the case {\bf (a)} and have convergence to an $\alpha$-stable distribution. 
On the other hand, in the LRD dominant case, $\frac1\alpha<2-2\sigma$, from Theorem \ref{our-auto}) we have almost sure convergence for $\frac1p>2-2\sigma$,  
yet for $\frac1p=(2-2\sigma)$ we have convergence to Rosenblatt process by {\bf (b)} .

Hence, Theorem \ref{our-auto} shows the $a.s$ convergence for difference of sample autocovariance and population autocovariance with  HT and LRD.
One example can be in the case that $h=0$. Theorem \ref{HTLRD} and (\ref{as-auto}) give the convergence in distribution
\begin{eqnarray}\nonumber
&& \frac{1}{a_n^2}\sum\limits_{k=1}^{n}(x_k^2-Ex_0^2)\stackrel{d}{\rightarrow} (S-\frac{\alpha}{\alpha-1})\sum\limits_{l=0}^{\infty}c_{l}^2\\\nonumber
&&\frac{1}{n^{2-2\sigma}}\sum\limits_{k=1}^{n}(x_k^2-Ex_0^2)\stackrel{d}{\rightarrow} \iota c_\sigma^2 U_\sigma(1),
\end{eqnarray}
for $\frac1p=\frac1\alpha$ and $\frac1p=2-2\sigma$, 
respectively. 

While, Theorem \ref{our-auto} gives the  almost sure convergence for  $\displaystyle\frac1{n^{\frac1p}}\sum\limits_{k=1}^{n}\left(x_k^2-Ex_0^2\right)$ when $\frac1p> (2-2\sigma)\vee \frac1\alpha$.\\
When we have convergence in distribution to non-trivial limits we can not get almost sure convergence to $0$. However, by Theorem \ref{our-auto} 
we can get arbitrary close to that with polynomial rate and get optimal polynomial almost sure rate of convergence. We can not do better than that in terms of MSLLN.

\section{Proofs} 

\subsection{{\bf Notation List}}
$|x|$ is Euclidean distance of some $\mathbb R^d$-vector $x$.\\
$\|C\|=\sup_{|x|=1}|Cx|$ for any $\mathbb R^{n\times m}$-matrix $C$.\\ 
$\lfloor t\rfloor\doteq \max\{i\in\mathbb N_0:i\le t\}$ and
$\lceil t\rceil\doteq \min\{i\in\mathbb N_0:i\ge t\}$ for any $t\ge 0$.\\
$a_{i,k}\stackrel{i}\ll b_{i,k}$ means that for each $k$ there is a $c_k>0$ that does not
depend upon $i$ such that $|a_{i,k}|\le c_k |b_{i,k}|$ for all $i,k$.\\
$\prod\limits_{l=p}^q B_l$ ($\forall\, B_l$ being a $R^{d\times d}$-matrix) $=B_qB_{q-1}\cdots B_p$ if $q\ge p$ or $I$ if $p>q$.\\
$a\vee b=\max\{a,b\}$ and
$a\wedge b=\min\{a,b\}$.

\subsection{ {\bf A First Light Tail Result} }

We first give a result that only handles long-range dependence without heavy tails.
However, our proof of Theorem \ref{HeavyLong} to follow will show that these 
two phenomena decouple, so we can easily build upon
the Theorem \ref{longonly} to handle both long-range dependence and heavy tails together.
\begin{theorem}\label{longonly}
Let $\left\{(\xi_{l},\overline \xi_{l}),\ {l\in\mathbb Z}\right\}$ be i.i.d.\ zero-mean 
random variables such that 
$E[(1+\xi_1^2)(1+\overline \xi_1^2)]<\infty$, $(c_{l},\overline c_{l})_{l\in\mathbb Z}$ satisfy
\begin{eqnarray}\nonumber
\sup\limits_{l\in\mathbb Z}|l|^\sigma|c_l|<\infty,\,\quad
\sup\limits_{l\in\mathbb Z}|l|^{\overline \sigma}|\overline c_l|<\infty \quad \mbox{for some}\quad \sigma, \overline \sigma \in \left(\frac{1}{2},1\right],
\end{eqnarray} 
\(
x_{k}=\sum\limits _{l=-\infty }^{\infty }c_{k-l}\xi _{l}\),
\(
\overline x_{k}=\sum\limits _{l=-\infty }^{\infty }\overline c_{k-l}\overline \xi _{l}\),
$d_k=x_k\overline x_k=\sum\limits_{l,m=-\infty}^{\infty}c_{k-l}\overline c_{k-m}\xi_{l}\overline \xi_{m}$ 
and $d=E[\xi_1\overline \xi_1]\sum\limits _{l=-\infty}^{\infty}c_{k-l}\overline c_{k-l}=E[\xi_1\overline \xi_1]\sum\limits_{l=-\infty }^{\infty }c_{l}\overline c_{l}$. 
Then, for $p<\frac1{2-\sigma-\overline\sigma}$
\[
\lim\limits_{n\rightarrow \infty}\frac{1}{n^{\frac{1}{p}}}\sum\limits_{k=1}^{n}\left(d_{k}-d\right)=0\ \ \mbox{a.s.}
\]
\end{theorem}

\proof.
Insomuch as the proof of the general case only differs cosmetically
from the notationally-simpler case where $\overline\xi_l=\xi_l$ and
$\overline c_l=c_l=\left\{ \begin{array}{ll} 1 & l=0\\|l|^{-\sigma}& l\neq 0\end{array}\right.$,
we only provide the proof of the later for which the constraint becomes $ p<\frac{1}{2-2\sigma}$. 
Assume without loss of generality that $\sigma <1$ and $E[\xi_1^2]=1$.

{\bf Step 1:} Divide partial sums into diagonal, large $c$, small and mixed type terms.\\ 
Let $n_{r}=2^{r}$ and $T=T\left(n\right)=n^{\nu}$ for $\nu>0$, $n\in\left[n_{r},n_{r+1}\right)$ 
and $r\in\mathbb{N}_{0}$, and define
\begin{eqnarray}
S^{\left(1\right)}_n\label{S1}
& = & \sum_{k=1}^n\sum\limits _{l=-\infty }^{\infty}c_{k-l}^2\left(\xi _{l}^{2}-1\right)
\\
S_n^{\left(2\right)}\label{S2}
& = & \sum_{k=1}^n 
\sum\limits_{
\stackrel{\scriptstyle l,m=k-T}{l\ne m}
}^{k+T}c_{k-l}c_{k-m}\xi _{l}\xi _{m}
\\S_{n}^{\left(3\right)}\label{S3}
& = & \sum_{k=1}^n 
\sum\limits _{\stackrel{\scriptstyle (l-k)\wedge(m-k)>T}{l\ne m}}c_{k-l}c_{k-m}\xi _{l}\xi _{m}\\
S_{n}^{\left(4\right)}\label{S4}
& = &\sum_{k=1}^n\sum\limits_{m-k>T}\sum\limits_{l=k-T}^{k+T}c_{k-l}c_{k-m}\xi _{l}\xi _{m}.
\end{eqnarray}
By breaking $\left\{\frac{1}{n^{\frac{1}{p}}}\sum\limits_{k=1}^{n}\left(d_{k}-d\right),\
n=1,2,...\right\}$
into pieces and considering those pieces with different (process) distributions,
we just need to show that
\[
\lim\limits_{n\rightarrow\infty}\frac{S_{n}^{\left(1\right)}}{n^{\frac{1}{p}}}
=\lim\limits_{n\rightarrow\infty}\frac{S_n^{\left(2\right)}}{n^{\frac{1}{p}}}
=\lim\limits_{n\rightarrow\infty}\frac{S_n^{\left(3\right)}}{n^{\frac{1}{p}}}
=\lim\limits_{n\rightarrow\infty}\frac{S_n^{\left(4\right)}}{n^{\frac{1}{p}}}
=0\ \mbox{a.s.,} 
\]
provided $ p<\frac{1}{2-2\sigma} $. 
To handle (the diagonal terms) $S_n^{(1)}$, we let $\zeta_l = \xi^2_l-1$, set $K=E[\zeta_1^2]$ and use standard steps.

{\bf Step 2:} Bound second moment of geometric diagonal partial sums $S^{(1)}_{n_r}$.\\
By symmetry and then integral approximation, we have that

\small{\begin{eqnarray}\label{Sn1}\nonumber
&& E[(S_{n_r}^{\left(1\right)})^2]\\\nonumber
& = & \nonumber
\sum\limits_{l=-\infty}^{\infty} \sum\limits_{m=-\infty}^{\infty}
\sum\limits _{j=1}^{n_r}\sum\limits _{k=1}^{n_r}
c_{k-l}^{2}c_{j-m}^{2}E[\zeta_{l}\zeta_{m}]
\\ 
&= &
K\sum\limits_{l=-\infty}^{\infty} 
\left|\sum\limits _{k=1}^{n_r} c_{k-l}^{2}\right|^2\nonumber
\\ \nonumber
&\stackrel{r}{\ll} &
\sum\limits _{k=1}^{n_r}\left(1+2\sum\limits_{l=1}^{\infty} l^{-4\sigma}
+2\sum\limits _{j=k+1}^{n_r}
\left(2(j-k)^{-2\sigma}+\sum\limits_{l=-\infty}^{k-1}(k-l)^{-2\sigma}(j-l)^{-2\sigma} 
\right.  \right.
\\ \nonumber
&+&\left.\left.\sum\limits_{l=k+1}^{j-1}(l-k)^{-2\sigma}(j-l)^{-2\sigma} 
+\sum\limits_{l=j+1}^{\infty}(l-k)^{-2\sigma}(l-j)^{-2\sigma}\right) \right) 
\\ \nonumber
&\stackrel{r}{\ll}&
\sum\limits _{k=1}^{n_r}\left(1+\sum\limits_{j=k+1}^{n_r}((j-k)^{-2\sigma}
+(j-k)^{1-4\sigma})\right)
\\ \,\,\quad
&\stackrel{r}{\ll}&
n_r.
\end{eqnarray}}
\normalsize
{\bf Note:}
\begin{eqnarray} \nonumber
\sum\limits _{l=k+1}^{j-1}\frac{1}{\left(l-k\right)^{2\sigma }\left(j-l\right)^{2\sigma }}
&\le& 
2\sum\limits _{l=k+1}^{\left\lfloor \frac{j+k}{2}\right\rfloor }\frac{1}{\left(l-k\right)^{2\sigma }\left(j-l\right)^{2\sigma }}\\\nonumber
&\stackrel{j,k}{\ll}& 
\left(j-k\right)^{-2\sigma }\sum\limits _{l=k+1}^{\left\lfloor \frac{j+k}{2}\right\rfloor }\frac{1}{\left(l-k\right)^{2\sigma }}\\
& \stackrel{j,k}{\ll}&
 \left(j-k\right)^{\left(1-4\sigma \right)}.
\end{eqnarray}

{\bf Step 3:} Maximal bound for geometric diagonal partial sums.\\
Following (\ref{Sn1}) we have for $n_{r}\le n<o<n_{r+1}$ 
\begin{eqnarray}\nonumber
E[(S_{o}^{\left(1\right)}-S_{n}^{\left(1\right)})^2]
&\le &
K\sum\limits_{l=-\infty}^{\infty} 
\left|\sum\limits _{k=n+1}^{o} c_{k-l}^{2}\right|^2\nonumber
\\ \nonumber
&\stackrel{o,n}{\ll}&
\sum\limits _{k=n+1}^{o}\left(1+\sum\limits_{j=k+1}^{o}((j-k)^{-2\sigma}
+(j-k)^{1-4\sigma})\right)
\\ 
&\stackrel{o,n}{\ll}&
o-n.
\end{eqnarray}
Therefore, it follows by Theorem 2.4.1 of Stout \cite{25} with $g(a,n)=Cn$ 
for some constant $C>0$ that
\begin{eqnarray}\label{max1}\nonumber
E\left[\max_{n_{r}\le n<o<n_{r+1}}\left(S_{o}^{(1)}-S_{n}^{(1)}\right)^{2}\right]
&\stackrel{r}{\ll}&\left(\frac{\log(2(n_{r+1}-n_r))}{\log2}\right)^2(n_{r+1}-n_{r})\\
&\stackrel{r}\ll& r^2n_r.
\end{eqnarray}

{\bf Step 4:} Use previous two steps to show \emph{normalized} diagonal sums converge.\\
Combining (\ref{Sn1}) and (\ref{max1}), one has that
\begin{eqnarray}\label{combine1}
\displaystyle \sum_{r=0}^\infty E\left[\max_{n_{r}\le n<n_{r+1}}\left(\frac{S_{n}^{(1)}}{n^\frac1p}\right)^{2}\right]
& \ll & \sum\limits_{r=0}^\infty r^2n_r^{1-\frac2p}<\infty,
\end{eqnarray}
provided ${p}\in(0,2)$.
It follows by Fubini's Theorem and $n^{\mbox{th}}$ term divergence that 
$$\lim\limits _{n\rightarrow\infty}\frac{S^{(1)}_n}{n^{\frac{1}{p}}}=0.$$

{\bf Step 5:} Set up for off-diagonal terms.\\ 
Letting 
\begin{eqnarray}
a^{2,n}_{l,m}&=&2\sum\limits _{k=1}^{n}1_{m-T\le k\le l+T}c_{k-l}c_{k-m}\\
a^{3,n}_{l,m}&=&2\sum\limits _{k=1}^{n}1_{k< l-T}c_{k-l}c_{k-m}\\
a^{4,n}_{l,m}&=&\sum\limits _{k=1}^{n}1_{k< m-T}1_{l-T\le k\le l+T}c_{k-l}c_{k-m},
\end{eqnarray}
we find that
\begin{eqnarray}
E\left[(S_{n}^{(i)})^{2}\right]\nonumber
& = & 
\sum\limits_{l_{1}=-\infty}^{\infty}\sum\limits_{m_{1}=l_{1}+1}^{\infty}a_{l_{1},m_{1}}^{i,n}
\sum\limits_{l_{2}=-\infty}^{\infty}\sum\limits_{m_{2}=l_{2}+1}^{\infty}a_{l_{2},m_{2}}^{i,n}
E\left[\xi _{l_{1}}\xi _{m_{1}}\xi _{l_{2}}\xi _{m_{2}}\right]
\\
& = & \nonumber
\sum\limits_{l_{1}=-\infty}^{\infty}\sum\limits_{m_{1}=l_{1}+1}^{\infty}a_{l_{1},m_{1}}^{i,n}
\sum\limits_{l_{2}=-\infty}^{\infty}\sum\limits_{m_{2}=l_{2}+1}^{\infty}a_{l_{2},m_{2}}^{i,n}
\delta_{l_1,l_{2}}
\delta_{m_1,m_{2}}
\\
& = & \sum\limits_{l=-\infty}^{\infty}\sum\limits_{m=l+1}^{\infty }\left(a_{l,m}^{i,n}\right)^{2}
\end{eqnarray}
and for $n_r\le n<o< n_{r+1}$
\begin{eqnarray}
E\left[(S_{o}^{(i)}-S_{n}^{(i)})^{2}\right]
=\sum\limits_{l=-\infty}^{\infty}\sum\limits_{m=l+1}^{\infty }\left(a_{l,m}^{i,o}-a_{l,m}^{i,n}\right)^{2}
\end{eqnarray}
for $i=2,3,4$.
Using a change of variables and the Beta distribution pdf, we have that
\begin{eqnarray}\nonumber
\sum\limits _{l=k+1}^{j-1}c_{j-l}c_{k-l}
&\stackrel{j,k}\ll&\label{betabound}
\int_{k}^j\left(j-t\right)^{-\sigma}\left(t-k\right)^{-\sigma}dt \\
&=&\left(j-k\right)^{1-2\sigma} 
\mathop{\underbrace{\int_{0}^1\left(1-s\right)^{-\sigma } s^{-\sigma }ds}}\limits_{B(1-\sigma,1-\sigma)}\stackrel{j,k}\ll \left(j-k\right)^{1-2\sigma}. 
\end{eqnarray}

{\bf Step 6:} Apply $S^{(1)}$-procedure for convergence of large $c$ terms $\frac{S_n^{(2)}}{n^\frac1p}$.\\ 
Using (\ref{betabound}) and integral approximation, one has for $n\in[n_r,n_{r+1})$
\begin{eqnarray}\nonumber
&& \!\!\!E\left[(S_{n}^{(2)})^2\right]
-4\sum\limits_{k=1}^{n}\sum\limits_{m>l}1_{k-T\le m\le k+T}\cdot 1_{k-T\le l\le k+T}c_{k-l}^{2}c_{k-m}^{2}
\\ \nonumber
& =& \!\!\! 8\sum\limits _{j>k}\sum\limits _{m>l}1_{j-T\le m\le k+T}\cdot 1_{j-T\le l\le k+T}c_{j-l}c_{j-m}c_{k-l}c_{k-m}
\\ \nonumber
 &  \le &\!\!\!  4\sum\limits _{j>k}\left|\sum\limits _{l=j-T}^{k+T}c_{j-l}c_{k-l}\right|^{2}
\\ \nonumber
 &  \le &\!\!\!  4\sum\limits _{k=1}^{n}\sum\limits _{j=k+1}^{n\wedge (k+2T)}
\left|2c_{j-k}+\sum\limits _{l=j-T}^{k-1}c_{j-l}c_{k-l}+\sum\limits _{l=k+1}^{j-1}c_{j-l}c_{k-l}+\sum\limits _{l=j+1}^{k+T}c_{j-l}c_{k-l}\right|^{2}
\\ \nonumber
 &  \stackrel{n}\ll  &\!\!\!  \sum\limits _{k=1}^{n}\sum\limits _{j=k+1}^{k+2T}\left[(j-k)^{-2\sigma}+\left(j-k\right)^{2-4\sigma }+\left(j-k\right)^{-2\sigma }T^{2-2\sigma }\right]
\\\nonumber
 & \!\!\! \stackrel{n}\ll&  nl(n),
\end{eqnarray}
where $l\left(n\right)
=\left\{ \begin{array}{ll}
T^{3-4\sigma}=n_r^{\nu(3-4\sigma)}&\sigma <\frac{3}{4}\\
\log \left(T\right)=\nu\log(n_r)&\sigma =\frac{3}{4}\\1&\sigma >\frac{3}{4}
\end{array}\right.$.
Hence,
\begin{eqnarray}\label{single2}
E\left[(S_{n}^{(2)})^2\right]
\stackrel{n}\ll  
nl\left(n\right) +\sum\limits_{k=1}^n\left|\sum\limits_{l=-T}^Tc_l^{2}\right|^{2}
\stackrel{n}\ll nl\left(n\right).
\end{eqnarray}
Similarly, we have for $ n_{r}\le n<o<n_{r+1}$ that
\begin{eqnarray}\nonumber
E\left[\left(S_{o}^{(2)}-S^{(2)}_{n}\right)^{2}\right]
& \stackrel{o,n}\ll &  
\sum\limits _{k=n+1}^{o}\left|\sum\limits _{l=-T}^{T}c_{l}^2\right|^{2}+
\sum\limits _{\stackrel{\scriptstyle j,k=n+1}{j> k}}^{o}\left|\sum\limits _{l=j-T}^{k+T}c_{j-l}c_{k-l}\right|^{2}
\\
&  \stackrel{o,n}\ll &  \left(o-n\right)l\left(n\right).
\end{eqnarray}
Therefore, it follows by Theorem 2.4.1 of Stout that
\begin{eqnarray}\label{Stout2}\nonumber
E\left[\max_{n_{r}\le n<o<n_{r+1}}(S^{(2)}_{o}-S^{(2)}_{n})^{2}\right]
&\stackrel{r}{\ll}&\left(\frac{\log(2n_r)}{\log2}\right)^2\!(n_{r+1}-n_{r})l(n_{r+1})\\
&\stackrel{r}\ll& r^2n_rl(n_{r}).
\end{eqnarray}
Combining (\ref{single2}) with $n=n_r$ and (\ref{Stout2}), one has that
\begin{eqnarray}
\displaystyle E\left[\sum_{r=0}^\infty\max_{n_{r}\le n<n_{r+1}}\left(\frac{S^{(2)}_{n}}{n^\frac1p}\right)^{2}\right]
& \ll & \sum\limits_{r=0}^\infty r^2n_r^{1-\frac2p}l(n_{r})<\infty,
\end{eqnarray}
provided $1+\nu(3-4\sigma)\vee0<\frac{2}{p}$ 
(i.e. $p<\frac{2}{1+\nu(3-4\sigma)}$ when $ \sigma<\frac34$ and $p<2$ when $\sigma\ge\frac34$, both of which are true).  
It follows that $\lim\limits_{n\rightarrow\infty}\frac{S^{(2)}_{n}}{n^{\frac{1}{p}}}=0$ a.s.

\vspace{1.5cm}
{\bf Step 7:} Apply $S^{(1)}$-procedure for convergence of small $c$ terms $\frac{S_n^{(3)}}{n^\frac1p}$.
\begin{eqnarray}\nonumber
\quad&\!\!\!\!\!\!&\!\! E\left[(S_n^{(3)})^{2}\right]\label{single3}\\ \nonumber
&  \!\!\!= \!\!\!& \!\!
8\sum\limits _{j>k}\sum\limits _{m>l}1_{j+T<l}\cdot 1_{k+T<l}c_{j-l}c_{j-m}c_{k-l}c_{k-m}\\\nonumber
&\!\!\!+\!\!\!&\!\! 4\sum\limits _{k=1}^{n}\sum\limits _{m>l}1_{k+T<l}c_{k-l}^{2}c_{k-m}^{2}
\\ 
& \!\!\!\le \!\!\!&\!\! \nonumber
4\sum\limits _{j>k}\left|\sum\limits _{l=j+T+1}^{\infty }c_{j-l}c_{k-l}\right|^{2}+2\sum\limits _{k=1}^{n}\left|\sum\limits _{l=k+T+1}^{\infty }c_{k-l}^{2}\right|^{2}
\\ 
& \!\!\!\stackrel{n}\ll \!\!\!&\!\!\nonumber
\sum\limits_{k=1}^{n-1}\sum\limits_{j=k+1}^{n}\left|\ \int_{j+T}^\infty\left(t-j\right)^{-\sigma }\left(t-k\right)^{-\sigma }dt \right|^{2}+\sum\limits _{k=1}^{n}\left|\ \int_{k+T}^\infty\left(t-k\right)^{-2\sigma }dt \right|^{2}
\\ 
 &  \!\!\!\stackrel{n}\ll \!\!\!&\!\! \nonumber
\sum\limits_{k=1}^n\left(\sum\limits_{j=k+1}^n\left|\ \int_{T}^\infty t^{-2\sigma}dt\right|^2
+\left|\ \int_{T}^\infty t^{-2\sigma}dt\right|^2\right)
\\  \quad
&\!\!\! \stackrel{n}\ll \!\!\!&\!\! n^2T^{2-4\sigma}.
\end{eqnarray}
Similarly, we have for $ n_{r}\le n<o<n_{r+1}$ that
\begin{eqnarray}
E\left[\left(S^{(3)}_{o}-S^{(3)}_{n}\right)^{2}\right]
\stackrel{o,n,r}\ll  \left(o-n\right)oT^{2-4\sigma}
\stackrel{o,n,r}\ll  \left(o-n\right)n_{r+1}^{1+\nu(2-4\sigma)}.
\end{eqnarray}
Therefore, it follows by Theorem 2.4.1 of Stout that
\begin{eqnarray} \label{Stout3}\nonumber
\!\!\!E\left[\max_{n_{r}\le n<o<n_{r+1}}\left(S_{o}^{(3)}-S_{n}^{(3)}\right)^{2}\right]
&\stackrel{r}\ll&\left(\frac{\log(2n_r)}{\log2}\right)^2
(n_{r+1}-n_{r})n_{r+1}^{1+\nu(2-4\sigma)}\\
&\stackrel{r}\ll& r^2 n_{r}^{2+\nu(2-4\sigma)}.
\end{eqnarray}
Combining (\ref{single3}) with $n=n_r$ and (\ref{Stout3}), one has 
\begin{eqnarray}
\displaystyle E\left[\sum_{r=0}^\infty\max_{n_{r}\le n<n_{r+1}}\left(\frac{S_{n}^{(3)}}{n^\frac1p}\right)^{2}\right]
& \ll &  \sum\limits_{r=0}^\infty r^2 n_r^{2+\nu(2-4\sigma)-\frac2p}<\infty,
\end{eqnarray} 
provided $p<\frac{1}{1+\nu(1-2\sigma)}$, which is the given condition, so $\lim\limits_{n\rightarrow\infty}\frac{S^{(3)}_{n}}{n^{\frac{1}{p}}}=0$ a.s.. 
\\

It is notable that condition on $p$,  $p<\frac{2}{1+\nu(3-4\sigma)}$, in step 6 gets more stringent when $\nu>1$ and the same 
is true for condition on $p$, $p<\frac{1}{1+\nu(1-2\sigma)}$, in step  7 when $\nu<1$, so the best choice that raises the same condition on $p$ is 
when $\nu=1$. Hence, we will have to satisfy $p<\frac{1}{1-2\sigma}$ in either cases.

{\bf Step 8:} Apply $S^{(1)}$-procedure for convergence of mixed terms $\frac{S_n^{(4)}}{n^\frac1p}$.\\ 
Finally, we note
\small{\begin{eqnarray}\nonumber
\!\!\!&\!\!\!\!\!\!&E\left[(S_{n}^{(4)})^2\right]\\\nonumber
\!\!\!&\!\!\!=\!\!\!&\sum\limits_{k=1}^{n}\sum\limits_{m=k+T+1}^\infty c_{k-m}^{2}\sum\limits_{l=k-T}^{l=k+T}c_{k-l}^{2}
 + 2\sum\limits _{k=1}^{n}\sum\limits _{j=k+1}^{k+2T}
\sum\limits_{m=j+T+1}^{\infty}c_{j-m}c_{k-m}
\sum\limits_{l=j-T}^{k+T}c_{j-l}c_{k-l}\\ \nonumber
\!\!\!& \!\!\!\stackrel{n}\ll\!\!\! &  
\sum\limits _{k=1}^{n}\left\{T^{1-2\sigma}+\sum\limits _{j=k+1}^{k+2T}
T^{1-2\sigma}
\left[(j-k)^{-\sigma}+\left(j-k\right)^{1-2\sigma }+\left(j-k\right)^{-\sigma }T^{1-\sigma }\right]
\right\}
\\\nonumber
\!\!\!&\!\!\! \stackrel{n}\ll\!\!\! &  nT^{3-4\sigma}.
\end{eqnarray}}
\normalsize
Similarly, we have for $ n_{r}\le n<o<n_{r+1}$ that
\begin{eqnarray}
E\left[\left(S_{o}^{(4)}-S^{(4)}_{n}\right)^{2}\right]
  \stackrel{o,n}\ll  \left(o-n\right)T^{3-4\sigma}\nonumber.
\end{eqnarray}
Therefore, it follows by $\nu=1$ and Theorem 2.4.1 of Stout that
\[
\!\!\!E\left[\max_{n_{r}\le n<o<n_{r+1}}\left(S^{(4)}_{o}-S^{(4)}_{n}\right)^{2}\right]
\stackrel{r}{\ll}\left(\frac{\log(2n_r)}{\log2}\right)^2(n_{r+1}-n_{r})n_{r+1}^{3-4\sigma}\stackrel{r}\ll r^2n_{r}^{4-4\sigma}.
\]
Combining these two equations, one has 
\begin{eqnarray}
\displaystyle E\left[\sum_{r=0}^\infty\max_{n_{r}\le n<n_{r+1}}\left(\frac{S^{(4)}_{n}}{n^\frac1p}\right)^{2}\right]
& \ll & \sum\limits_{r=0}^\infty r^2n_r^{(4-4\sigma)-\frac2p}<\infty,
\end{eqnarray}
provided $p<\frac{1}{2-2\sigma}$, which is true.
It follows that $\lim\limits_{n\rightarrow\infty}\frac{S^{(4)}_{n}}{n^{\frac{1}{p}}}=0$ a.s.
\eproof
\subsection{{\bf Proof of Theorem \ref{HeavyLong}}}
Without loss of generality we assume $1<\alpha<2$.

{\bf Step 1:} Reduce to continuous $\{(\xi_l,\overline \xi_l)\}$.\\
Let $\{(U_l)\}_{l\in\mathbb Z}$ be independent 
$[-1,1]$-uniform random variables that are independent of everything
and set $\overline U_l=U_l$ for all $l$.
Then, we have that
\small
\begin{eqnarray}\nonumber
\frac{1}{n^{\frac{1}{p}}}\sum\limits_{k=1}^{n}\left(d_{k}-d\right)&=& 
\frac{1}{n^{\frac{1}{p}}}\sum\limits_{k=1}^{n}
\sum\limits_{l,m=-\infty}^{\infty}c_{k-l}\overline c_{k-m}
\left((\xi_{l}+U_l)(\overline \xi_{m}+\overline U_m) -d-\frac23\right)\\
&-&\frac{1}{n^{\frac{1}{p}}}\sum\limits_{k=1}^{n}
\sum\limits_{l,m=-\infty}^{\infty}c_{k-l}\overline c_{k-m}
\left(\xi_{l}\overline U_m+U_l\overline \xi_{m}+U_l\overline U_m-\frac23\right).
\end{eqnarray}
\normalsize
However, 
\begin{eqnarray}
\lim\limits_{n\rightarrow \infty}
\frac{1}{n^{\frac{1}{p}}}\sum\limits_{k=1}^{n}
\sum\limits_{l,m=-\infty}^{\infty}c_{k-l}\overline c_{k-m}
\left(\xi_{l}\overline U_m+U_l\overline \xi_{m}+U_l\overline U_m-\frac23\right)=0
\end{eqnarray}
by Theorem \ref{longonly}.
Moreover, $\xi_1+U_1,\overline\xi_1+\overline U_1$ have the same moment and tail probability
bounds as $\xi_1,\overline\xi_1$.
Hence, without loss of generality, we can assume $\xi_l,\overline \xi_m$ are continuous
random variables, which will be important for the truncation to follow in Step 4.

{\bf Step 2:} Handle off-diagonal sum as previous proof since unaffected by heavy tails.\\
Suppose $S_n^{\left(2\right)}$, $S_n^{\left(3\right)}$ and $S_n^{\left(4\right)}$
are defined as in (\ref{S2}-\ref{S4}).
Then, we know that
\[
\lim\limits_{n\rightarrow\infty}\frac{S_n^{\left(2\right)}}{n^{\frac{1}{p}}}
=\lim\limits_{n\rightarrow\infty}\frac{S_n^{\left(3\right)}}{n^{\frac{1}{p}}}
=\lim\limits_{n\rightarrow\infty}\frac{S_n^{\left(4\right)}}{n^{\frac{1}{p}}}
=0\ \mbox{a.s.,} 
\]
provided $p<\frac{1}{2-\sigma -\overline \sigma}$ by the proof of 
Theorem \ref{longonly}. 

{\bf Step 3:} Reduce $\xi_l\overline \xi_l$ (in diagonal sum) to non-negative with single atom at $0$.\\
Noting 
\small{\begin{eqnarray}\label{posnegdecomp}\nonumber
&&\!\!\!\!\sum\limits_{l=-\infty}^\infty c_{k-l}\overline c_{k-l}(\xi_l\overline \xi_l -E[\xi_l\overline \xi_l ])
\\\quad
&=&\!\!\!\!\sum\limits_{l=-\infty}^\infty c_{k-l}\overline c_{k-l}((\xi_l\overline \xi_l)^+-E[(\xi_l\overline \xi_l)^+])
-\sum\limits_{l=-\infty}^\infty c_{k-l}\overline c_{k-l}((\xi_l\overline \xi_l)^--E[(\xi_l\overline \xi_l)^-]),
\end{eqnarray}}
\normalsize

\noindent we only have to consider the case where $\xi_l\overline \xi_l\ge 0$ for the remainder of the proof.
Moreover, insomuch as the proof of the general case only differs cosmetically
from the notationally-simpler case where $\overline\xi_l=\xi_l$, $E[\xi_1^2]=1$ and
$\overline c_l=c_l=\left\{ \begin{array}{ll} 1 & l=0\\|l|^{-\sigma}& l\neq 0\end{array}\right.$,
we only provide the proof of the later for which the long-range dependence
constraint becomes $p<\frac{1}{2-2\sigma}$. 
We will however indicate the most significant changes that would be needed for the general case.

{\bf Step 4:} Divide diagonal terms into zero-mean truncated (i.e.\ bounded) and remainder pieces.\\ 
Let $\kappa>0$. 
Fix $u_r^+= n_r^\frac{\kappa}{2-\alpha}$ to find
\begin{eqnarray}\label{uplus}
2\int_0^{u_r^+}P(\xi_1^2>s)sds\stackrel{r}{\ll}2\int_0^{u_r^+}s s^{-\alpha}ds
\stackrel{r}{\ll} n^\kappa_r\ \ \forall\ r=1,2,...
\end{eqnarray} 
Now, by defining
\begin{eqnarray}
\quad\left\{ \begin{array}{l} \overline \zeta_i=\overline \zeta_i^r=(\xi_i^2\wedge u_r^+)-\vartheta_i, \mbox{ where } 
\vartheta_i\doteq\int_0^{u_r^+}P(\xi_i^2>s)ds\leq 1, \\
\tilde{\zeta}_i= \tilde{\zeta}_i^r=\xi_i^2-1-\overline \zeta_i^r,
\end{array}\right.
\end{eqnarray}
we find that
\begin{eqnarray}
E[\overline \zeta_i]=\int_0^{u_r^+}P(\xi_i^2>t)dt-\int^{u_r^+}_0 P(\xi_i^2>t)dt=0,
\end{eqnarray}
so both $\overline \zeta_i$ and $\tilde \zeta_i$ are zero mean,
and by (\ref{uplus})
\begin{eqnarray}\label{barzerabound}\nonumber
E[|\overline\zeta_1|^2]\!\!&=&\!\!E|\xi^2_1\wedge u_r^+|^2-\left(\int^{u_r^+}_0 P(\xi^2_1>t)dt\right)^2
\\
\!\!&=&\!\! 2\int_0^{u_r^+}P(\xi^2_1>s)sds-\left(\int^{u_r^+}_0 P(\xi^2_1>t)dt\right)^2
\stackrel{r}{\ll} n^\kappa_r\ \ \forall\ r=1,2,...
\end{eqnarray}
(In the general case, we note that $\xi_1\overline \xi_1$ is non-negative and of continuous distribution
on $(0,\infty)$ so $E[\xi_1\overline \xi_1\wedge u_r^+]=\int_0^{u_r^+}P(\xi_1\overline \xi_1>s)ds$ as required.
We also have $\tilde{\zeta}_i^r=\xi_i\overline\xi_i-E[\xi_i\overline\xi_i]-\overline \zeta_i^r$.)

{\bf Step 5:} Moment Bound for truncated using the proof of Theorem \ref{longonly}.\\
Noting $\{\overline \zeta_i\}$ are i.i.d. with $E[\overline \zeta_1]=0$ and $E[\overline \zeta_1^2]<\infty$ and
defining
\begin{eqnarray}
S^{\left(1\right)}_n
&=&\sum_{k=1}^n\sum\limits_{l=-\infty }^{\infty}c_{k-l}^2\overline\zeta_{l}, 
\end{eqnarray}
one finds 
from (\ref{combine1}) in the proof of Theorem \ref{longonly} that
\begin{eqnarray}\label{max1bar}
E\left[\max_{n_{r}\le n<n_{r+1}}\left(S_{n}^{(1)}\right)^{2}\right]
\le E|\overline \zeta_1|^2 r^2n_r.
\end{eqnarray}
Hence, it follows by (\ref{barzerabound}) that
\begin{eqnarray}\label{max1bar2}
E\left[\max_{n_{r}\le n<n_{r+1}}\left(S_{n}^{(1)}\right)^{2}\right]
\stackrel{r}{\ll} 
r^2 n^{1+\kappa}_r.
\end{eqnarray}

{\bf Step 6:} Moment Bound for remainder using Doob's inequality.\\
Turning to the $\tilde \zeta_i^r$
and using the formula
\begin{eqnarray}\label{expform}
E[g(X)]=\int_0^\infty g'(t)P(X>t)dt-\int_{-\infty}^0 g'(t)P(X<t)dt,
\end{eqnarray}
one has by our tail probability bounds that the non-negative part of $\tilde \xi_1$
satisfies
\begin{eqnarray}\label{momboundtilde}\nonumber
E|\tilde \zeta_1^+|^\tau
&=&
\tau\int_0^\infty s^{\tau-1}P(\xi_1^2>u_r^++s+1-\vartheta_1)ds\\\nonumber
&\leq&
\tau\int_0^\infty s^{\tau-1}P(\xi_1^2>u_r^++s)ds\ \mbox{since }\vartheta_1\le1\\\nonumber
&\stackrel{r}{\ll}& 
\int_{u_r^+}^\infty (s-u_r^+)^{\tau-1}s^{-{\alpha}}ds\\\nonumber
&\le& \int_{u_r^+}^{2u_r^+}(s-u_r^+)^{\tau-1}ds (u_r^+)^{-{\alpha}}
+ \int_{2u_r^+}^\infty (s-u_r^+)^{\tau-\alpha-1}ds\\
&\stackrel{r}{\ll}& (u_r^+)^{\tau-{\alpha}}
\stackrel{r}{\ll} n_r^{\frac{\kappa(\tau-{\alpha})}{2-\alpha}},
\end{eqnarray}
for $1<\tau<\alpha$.
Therefore, it follows by Jensen's inequality and Doob's $L_p$ inequality that
\begin{eqnarray}\nonumber
E^\frac1\tau\left[\sup\limits_{n_r\le n< n_{r+1}}\left|\sum_{k=1}^n
\sum_{l=-\infty}^\infty c_l^2\tilde \zeta^r_{k-l}\right|^\tau\right]
\!\!\!&\le&\!\!\!
E^\frac1\tau\left[\left|
\sum_{l=-\infty}^\infty c_l^2\sup\limits_{n_r\le n< n_{r+1}}\left|\sum_{k=1}^n
\tilde \zeta_{k-l}^r\right|\right|^\tau\right]\label{alphaterm}\\
\!\!\!&\stackrel{r}{\ll}&\!\!\!
\sum_{l=-\infty}^\infty c_l^2
E^\frac1\tau\left[
\sup\limits_{n_r\le n< n_{r+1}}\left|\sum_{k=1}^n
\tilde \zeta_{k-l}^r\right|^\tau\right]\nonumber\\
\!\!\!&\stackrel{r}{\ll}&\!\!\!
\sum_{l=-\infty}^\infty c_l^2
E^\frac1\tau\left[
\left|\sum_{k=1}^{n_{r+1}-1}
\tilde \zeta_{k-l}^r\right|^\tau\right]\nonumber\\
\!\!\!&\stackrel{r}{\ll}&\!\!\! n_r\|\tilde \zeta_1^r\|_{\tau},
\end{eqnarray}
so by (\ref{momboundtilde},\ref{alphaterm})
\begin{eqnarray}\label{remmoment}
E\left[\sup\limits_{n_r\le n< n_{r+1}}
\left|\sum_{k=1}^n
\sum_{l=-\infty}^\infty c_l^2\tilde \zeta^r_{k-l}\right|^\tau\right]
\stackrel{r}{\ll}
n_r^{\tau - \frac{\kappa(\alpha-\tau)}{2-\alpha}}.
\end{eqnarray}

{\bf Step 7:} Use Truncation and Error Term bounds with Borel-Cantelli for convergence.\\
Combining (\ref{max1bar2}) and (\ref{remmoment}), one has that
\small{\begin{eqnarray}\label{Pmaxzeta}\nonumber
&& P\left(\sup\limits_{n_r\le n< n_{r+1}}\left|
\sum_{k=1}^n
\sum_{l=-\infty}^\infty c_l^2\zeta_{k-l}\right|>2\epsilon n_r^\frac{1}p\right)
\\ &\le&
\frac{E\left[\sup\limits_{n_{r}\le n<n_{r+1}}
\left|\sum\limits_{k=1}^n
\sum\limits_{l=-\infty}^\infty c_l^2\overline \zeta^r_{k-l}\right|^{2}\right]}{\epsilon^2 n_r^\frac2p}
+\frac{E\left[\sup\limits_{n_{r}\le n<n_{r+1}}
\left|\sum\limits_{k=1}^n
\sum\limits_{l=-\infty}^\infty c_l^2\tilde \zeta^r_{k-l}\right|^\tau
\right]}{\epsilon^\tau n_r^\frac{\tau}p}\nonumber\\
&\stackrel{r}{\ll}&
r^2 n_r^{1+\kappa-\frac2p}+n_r^{\tau - \frac{\kappa(\alpha-\tau)}{2-\alpha}-\frac{\tau}p}\nonumber\\
\quad\quad&\!\!\!\stackrel{r}{\ll}&
r^2 n_r^{1-\frac{\alpha}p}+n_r^{\tau -\frac{\alpha}p} ,
\end{eqnarray}}
\normalsize
by letting $\kappa = \frac{2-\alpha}p$.
Hence, if $\tau\in\left(1,\frac{\alpha}p\right)$, then 
\begin{eqnarray}\label{Psumzeta}
\sum_{r=1}^\infty P\left(\sup\limits_{n_r\le n< n_{r+1}}\left|
\sum_{k=1}^n
\sum_{l=-\infty}^\infty c_l^2\zeta_{k-l}\right|>2\epsilon n_r^\frac{1}p\right)<\infty ,
\end{eqnarray}
under our heavy-tail condition $p<\alpha$ and
\begin{eqnarray}\label{aszeta}
n^{-\frac{1}p}
\sum\limits_{k=1}^n
\sum\limits_{l=-\infty}^\infty c_l^2\zeta_{k-l}\to 0\ \ \mbox{a.s.} ,
\end{eqnarray}
by Borel-Cantelli.
The proof is complete. \eproof

\end{document}